\documentclass[10pt]{article}
 \usepackage{amssymb}

\def\bC{{\mathbb C}}

\def\bQ{{\mathbb Q}}

\def\t{\tilde}

\def\E{\bf E}

\def\P{{\mathbb P}^1}  
\def\M{{\mathcal M}(Y)}

\newtheorem {Theorem}{Theorem}[section]

\newtheorem {Definition}[Theorem]{Definition}

\newtheorem {Remark}[Theorem]{Remark}

\newcommand{\beq}{\begin{equation}}
\newcommand{\eeq}{\end{equation}}

 \def\bC{{\mathbb C}}  \def\t{{\tau}}

    \def\g{{\sigma }}
 \def\M{{\mathcal M}}   \def\s{{\sigma }}

\def\Aut{\mbox {Aut}}
 \def\BB{{\cal B}}
\def\Bbr{{\cal B}_{r}}
\def\BBR{{\cal B}^{(r)}}

 \def\E{{\frak E}}

\begin{document}
\setcounter{section}{-1}

 \centerline{\Large A  GAP package for braid orbit computation, and applications}
         \bigskip

          \bigskip

         \centerline{K. Magaard
          \footnote{Partially supported by NSA grant MDA-9049810020},
          S. Shpectorov,  and Helmut V\"olklein
         \footnote{Partially supported by NSF grant DMS-0200225 }}

 \bigskip

  \centerline{Wayne State University, 
 Bowling Green State University and University
of Florida }

         \vspace{2cm}
        \noindent {\bf Abstract:} \ Let $G$ be a finite group.
By Riemann's Existence Theorem, braid orbits of generating systems 
of $G$ with product 1 correspond to irreducible
families of covers of the Riemann sphere with monodromy group $G$.
Thus many problems on algebraic curves require the computation of 
braid orbits. In this paper we describe an implementation of this
computation. We discuss several applications, including
the classification of  
 irreducible
families of indecomposable rational functions with exceptional
monodromy group.

\section{Introduction}

Let $G$ be a finite group and $\s=(\s_1,...,\s_r)$ a tuple of elements
of $G$ with $\g_1\cdots \g_r=1$.
The  {\bf braid orbit} of $\s$ 
is the smallest set of
tuples from $G$ that contains $\s$ and is closed under
the braid operations
$$ (g_1,...,g_r)^{Q_i}\ \ = \ \
(g_1,\ldots,\ g_{i+1},\ g_{i+1}^{-1}
 g_i g_{i+1} \ ,\ldots, g_r)\leqno{(1)} $$
for $i=1,...,r-1$. 
Clearly,
the unordered collection of conjugacy
classes $C_1,...,C_r$ represented by the elements of the tuple
is an invariant of the braid orbit.
This paper describes a package of programs
written in \cite{GAP4} for the computation of all braid orbits 
associated with
given classes $C_1,...,C_r$. We call it the BRAID program.
It is available at http://www.math.ufl.edu/~helmut. An alternative
approach has recently been worked out by Kl\"uners (Kassel),
 using MAGMA.
A precursor was the HO-program of Przywara \cite{Pr}
 which is now outdated.

Our interest in computing braid orbits comes from the fact that they
correspond to irreducible
families of covers of the Riemann sphere. This is a classical fact,
used by  Hurwitz (who found formula (1)) and  many algebraic
geometers since then. This connection to geometry is briefly explained 
in section
\ref{expl}.
 The version required  for the application
to the Inverse Galois Problem was worked out by Fried and V\"olklein
\cite{FV}. 

However, there are also
purely group-theoretic applications of our braid program, 
e.g., to find generators of
a given group with prescribed element orders.
Most applications have been in geometry and number theory, though,
via the connection to covers. Covers of $\P$
defined over $\bQ$ yield Galois
realizations of $G$ over $\bQ$ via Hilbert's irreducibility theorem
--- the braid program is needed to find suitable covers for which the
criteria of Inverse Galois Theory apply. A good example of that is
Malle's construction \cite{Ma}
 of multi-parameter polynomials with various  small Galois groups.
His $L_3(2)$-polynomial is used as an example in section \ref{L32}
below to obtain a generic rational function of degree 7 with
monodromy group $L_3(2)$. Another example is
Matzat's realization \cite{MM},III, 7.5, of $M_{24}$, for which  Granboulan
\cite{Gra} computed an explicit polynomial.  
A further example is the realization of symplectic groups $Sp(n,q)$
by Thompson and V\"olklein \cite{ThV} which depends on the fact that the
pure braid operations (2) generate an abelian group of permutations of
the corresponding braid orbit (mod conjugation).
There are numerous other  applications to the Inverse Galois Problem,
see \cite{MM} and \cite{Buch}.

There are also applications to problems about the geometry
of algebraic curves and their moduli spaces $\M_g$. E.g., 
in \cite{kyoto} the authors study the  locus in $\M_g$ of curves 
with given \lq large\rq\
automorphism group $G$. The irreducible 
components of that locus correspond
to certain  braid orbits in $G$. The BRAID program enabled us
to completely classify these components for $g\le10$
 and compute the genus
of those that are 1-dimensional.

In this paper we describe the application to classifying the irreducible
families of indecomposable rational functions with 
monodromy group other than $S_n$ or $A_n$. 
A generating system $\g_1,\ldots, \g_r$ of a 
transitive permutation group $G$
with $\g_1\cdots \g_r=1$ is called a {\bf genus zero system}
if the corresponding covers of $\P$ have genus $0$, i.e., are given 
by a rational function $f(x)\in\bC(x)$. The function $f$ is
{\bf indecomposable} (with respect to composition) if and only if
$G$ is primitive. In this case we say $\g_1,\ldots, \g_r$ is a
{\bf primitive genus zero system}.  There is a huge variety
of such systems that generate $S_n$ or $A_n$, too many
to be classified. 
 Those functions with smaller monodromy group
satisfy interesting identities and therefore it seems desirable to
have a complete classification of their irreducible families.

 Thus we need to compute all
braid orbits of genus zero systems in primitive
permutation groups $G$ other than
$A_n$ or $S_n$. 
It follows from the proof of the Guralnick-Thompson Conjecture
(see \cite{FM}) that only finitely many  groups $G$ occur.
The complete list is being worked out by 
Frohardt, Guralnick, Magaard and Shareshian 
\cite{FGM2}, \cite{GS} (project nearly completed).
 The smallest group that occurs
is $G=L_3(2)$ (acting on 7 points). We  study this example in
section \ref{L32}.
In section \ref{five} we present all braid orbits of
genus zero systems of length $\ge5$ in almost simple groups
other than
$A_n$ or $S_n$.
 The remaining cases
(length $3$ and $4$) will be collected in a data base, there is too
many of them to be displayed here.

Another application of the BRAID program was given in \cite{MV}.
  We say a tuple
$\g_1,\ldots, \g_r$ in $S_n$ has {\bf full moduli dimension} if the
corresponding family of covers contains the general curve of that
genus. If that holds and the genus is at least 4 then $\g_1,\ldots, \g_r$
generate $S_n$ or $A_n$ by work of Guralnick and others \cite{GM},
\cite{GS}. In genus 2 and 3 there are several other possible cases.
In  \cite{MV} it was shown that  the general curve of genus 3
has a cover to $\P$ of degree 7 with monodromy group $L_3(2)$.
  The associated  tuple consists  of 9 involutions (with product 1)
generating $L_3(2)$.  
 There is only one 
braid orbit of such tuples by \cite{MV}, Remark 5.1.
 This requires an iterative
application of the BRAID program because the  orbit
 is too large for a direct computation.  This iterative procedure
for computing  braid-orbits of long tuples in small groups requires
computing braid-orbits of (shorter) tuples of product $\ne1$ (see 
Remark \ref{rem}).


\part{Description of the BRAID program}

\section {Exact formulation of the problem}

Fix an integer $r\ge3$.

The Artin braid group $\Bbr$ is defined by a presentation
on generators $Q_1,...,Q_{r-1}$ and relations
$$Q_iQ_{i+1}Q_i\ = \ Q_{i+1}Q_iQ_{i+1}
\ \ \ \mbox{and} \ \ \
Q_iQ_j\ = \ Q_jQ_i \ \ \mbox{for} \ |i-j|>1$$
Mapping $Q_i$ to the transposition $(i,i+1)$ extends to a
homomorphism $\kappa: \Bbr\to S_r$
with kernel $\BBR$, the pure Artin braid group. It is generated by the
$$ Q_{ij} \ = \ Q_{j-1}\cdots Q_{i+1}\ Q_{i}^2\ Q_{i+1}^{-1}\cdots
Q_{j-1}^{-1} \ = \ Q_{i}^{-1}\cdots Q_{j-2}^{-1}\ Q_{j-1}^2\ Q_{j-2}
\cdots Q_i, \ \ \ 0\le i<j\le r\leqno{(2)} $$
More generally, if $P$ is a partition of $\{1,\dots,r\}$, let
$S_P$ be the stabilizer of $P$ in $S_r$ and set 
$\BB_P = \kappa^{-1}(S_P) $. We always choose $P$ such that each block
consists of all integers between the smallest and largest element of
the block. Thus we can identify $P$ with the list
 of the lengths of its parts. 
$\BB_P$ is generated by the $ Q_{ij}$ with $i,j$ not in the same
block of $P$, and the $ Q_{i}$ with $i,i+1$ in the same block.

Now let $G$ be a finite group.
Then $\Bbr$ acts on $r$-tuples  of elements of $G$ with product 1
via formula (1) above.
The orbits of this $\Bbr$-action are called {\bf braid orbits}.
This $\Bbr$-action commutes with the action of $\Aut(G)$ on tuples
defined by
$$\alpha(\g_1,\dots,\g_r) \ \ = \ \ (\alpha(\g_1),\dots,\alpha(\g_r))$$
 for $\alpha\in \Aut(G)$. Thus $\Bbr$ permutes $\Aut(G)$-orbits
(as well as Inn$(G)$-orbits) of tuples.

Note that in the $\Bbr$-action on tuples $(\g_1,\ldots,\g_r)$,
the conjugacy classes $\g_1^G,...,\g_r^G$ are being permuted
via the map $\kappa: \Bbr\to S_r$. This yields an obvious simplification
in computing the braid orbit of a tuple $(\s_1,\ldots,\s_r)$: We only need
to compute those tuples in the braid orbit where the
classes $\s_1^G,...,\s_r^G$ occur in that given order. In other words,
we only compute the orbit of $(\s_1,\ldots,\s_r)$ under the subgroup
of $\Bbr$ that stabilizes this order of the conjugacy classes. This subgroup
equals $\BB_P$, where $P$ is the partition of $\{1,\dots,r\}$ such that
$i$ and $j$ lie in the same block iff $\s_i$ is conjugate
$\s_j$. 

The classes $\s_1^G,...,\s_r^G$ have an important interpretation in terms
of the associated covers ("distinguished inertia group generators", see
\cite{Buch}). 
Thus we consider the
following basic problem.

\medskip\noindent
{\bf Problem 1:} {\sl Let $C_1,...,C_r$ be non-trivial
conjugacy classes of the
finite group $G$.
Let $P$ be the  partition of $\{1,\dots,r\}$ such that
$i$ and $j$ lie in the same block iff $C_i=C_j$. We want to
compute the orbits of $\BB_P$
on the set of Inn$(G)$-orbits on}
 $$\E(C_1,...,C_r) \ \ = \ \ \{(\g_1,\ldots,\g_r):\ \g_i\in C_i,\
\g_1\cdots \g_r=1\}$$

\medskip\noindent
Further geometric information is furnished by the permutations induced
by certain of the generators of $\BB_P$ on the braid orbit. So we record
these permutations as we construct the braid orbit. In the case $r=4$,
for example, this information can be used to compute the genus
of the corresponding Hurwitz curve
(see section \ref{genus} below).

\begin{Remark} \label{rem}
Modified versions of Problem 1 arise where $\BB_P$ is replaced
by a subgroup $\BB'$. For example, $\BB'$ could be $\BB_{P'}$ for a
partition $P'$ finer than $P$, or it could be an analogous
 subgroup of $\BB_{r-1}$. The latter is equivalent to acting on tuples of length
$r-1$ with product $\ne1$. (Note that the braid group  acts on
tuples with any fixed product by formula (1)). Further choices for $\BB'$
are the subgroups of the braid group induced by the fundamental groups of
certain curves on the configuration space, see \cite{Dett}; generators
for some of these groups can be found at\\ 
http://www.iwr.uni-heidelberg.de/groups/compalg/dettweil/papers.html.
(They have applications to the Inverse Galois Problem). The BRAID program
can easily be adapted to these modified versions of Problem 1.
\end{Remark}

\def\t{\tau}
\section {Program input and output} 
Problem 1 is solved by our main routine {\tt AllBraidOrbits}. To call this
routine, choose a tuple $\t$ representing the classes 
$C_1,...,C_r$. (The tuple $\t$ need not have product 1). The classes
$C_1,...,C_r$ must be ordered such that if $C_i=C_j$ with $i<j$ then
$C_i=C_k$ for all $i\le k\le j$.
The cardinality $c$ of $\E(C_1,...,C_r)$ is given
by a well-known formula (see \cite{MM}, Ch.I, Th. 5.8)
involving the values on $C_1,...,C_r$
of the irreducible characters of $G$. This number $c$ is called the
{\bf structure constant} associated with $C_1,...,C_r$. It can be 
computed with the GAP command {\tt ClassStructureCharTable}, once the
character table of $G$ is available.
  Once $c$ has been computed, we 
call our main routine in the form 
\begin{center}{\tt AllBraidOrbits("ProjectName",$G,\t,P,c$)} \end{center}
where {\tt ProjectName} is any string that is used to label the output
files. Here $G$ has to be a permutation group because
many standard algorithms of GAP4 work only in that case. The routine
computes the ${\BB_{P}}$-orbits on 
$\E(C_1,...,C_r)$ mod Inn$(G)$. 
 For each orbit it creates a file containing   a list of
representatives of Inn$(G)$-orbits of the tuples  in the 
orbit, plus the permutations induced on the orbit by the
generators of ${\BB_{P}}$ and by the generators of the pure braid group. 

\smallskip\noindent
{\bf User-friendly version:}  $G$ and $\tau$ are as above.
The routine
\begin{center}{\tt Braid($G$,$\t$)} \end{center}
firstly computes the character table of $G$ and uses it to compute
the structure constant  $c$. For large $G$ this may be time-consuming
or not feasible at all (then the character table must be taken from
some library). Furthermore, the program computes the partition $P$.
Then it calls {\tt AllBraidOrbits}, using always the same ProjectName
"TEMP". The previous contents of that directory is removed each time the
routine is called. In the end,
it summarizes the output by listing all braid orbits found that
consist of tuples $\s$ generating $G$. If $r=4$, the genus of the
inner Hurwitz
curve ${\cal H}_{\mbox{in}}^{\mbox{red}}(\s)$
and straight inner Hurwitz
curve $\tilde{\cal H}_{\mbox{in}}^{\mbox{red}}(\s)$
are given for each of those orbits (see section \ref{expl},  \ref{genus}).
A variation is the command
\begin{center}{\tt Braid($G$,$\t$,$U$)} \end{center}
where $U$ is a 
core-free subgroup of $G$ of index $n$.
Now the routine calls {\tt AllBraidOrbits} with $G$ replaced by its
 normalizer in $S_n$, where $G$ is embedded in $S_n$ via its
permutation representation on the cosets of $U$. If $r=4$, 
the genus of the Hurwitz
curve ${\cal H}^{\mbox{red}}(\s)$ (relative to this
permutation representation) is given for each orbit
of tuples generating $G$.

\section{Description of the algorithm}

At the beginning of its main loop, the {\tt AllBraidOrbits} routine
collects a batch of 
random tuples from $\E(C_1,\ldots,C_r)$.  If one of these tuples 
does not belong to a known (braid) orbit, 
a routine {\tt BraidOrbit} is called to 
generate the new orbit and add it to the 
list of known orbits. Furthermore,   
the variable $c$ is adjusted to be the number of tuples in 
$\E(C_1,\ldots,C_r)$ which do not belong to any one of the
currently known orbits. When $c=0$, we are done.

One is  mainly interested in those
tuples from $\E(C_1,\ldots,C_r)$ that
generate $G$. However, we don't know how to determine their
number beforehand (in any efficient way). That's why we are working
with the larger set $\E(C_1,\ldots,C_r)$ (whose cardinality $c$ is given
by the structure constant formula). Here are some variations
on choosing  the input value
of $c$: Setting  $c$
to a very large number, {\tt AllBraidOrbits} is turned into an
infinite loop. The user breaks the loop when he is convinced that all 
relevant orbits have been found. This avoids the actual computation
of the structure constant. On the other hand, by setting $c$
 below the actual size of $\E(C_1,\ldots,C_r;1)$ one can 
skip the last few small orbits that are usually irrelevant.  
For example, if only the orbits of 
generating tuples are of interest then one can  quit once 
the number of tuples unaccounted for is below $|G/Z(G)|$.

Hitting a particular small orbit with a random tuple is 
not likely to happen quickly. Therefore,
 we implemented 
a particular way of creating random tuples. It
involves maintaining a list of small subgroups  
generated by known tuples, and trying to find
more tuples in those subgroups. For example,
the case of $6$-tuples of
double transpositions
in $A_7$ took about 2 hours using a purely random tuple selection. Our
current method cut this time to
30 minutes. In both cases the
 program took 20 minutes to account for about $90
\%$ of the tuples.
So the time for finding the last $10 \%$  was cut
from 100 minutes to
10 minutes.

The routine {\tt BraidOrbit}($\s$) constructs the braid orbit
of a tuple $\s$.
 We use a Dixon-Schreier
algorithm: Beginning with $\s$, apply the  generators of ${\BB_{P}}$ one by one 
to the known tuples  and check wether or not the image 
is $G$-conjugate to one them. If not we append the new tuple to the list.
The routine terminates when  
 no further tuples can be produced.

The only difficulty is how to check efficiently whether two given tuples
are $G$-conjugate.
To speed this up we use a fingerprinting 
technique. Fingerprints  are sequences of numbers 
that can be quickly computed for a tuple. Tuples with distinct 
fingerprints cannot be conjugate. Currently, fingerprints are realized 
as the orders (as group elements) of certain random words in 
$\g_1,\ldots,\g_r$. 
The fingerprints are stored along with the tuples. 
Access to a tuple is via its fingerprint. Access to a fingerprint is 
via a hash table, the address for which is formed from the entries of 
the fingerprint. We remark that this method works well for a large 
variety of groups $G$. Exceptions are  
Frobenius groups and some $p$-groups.

\section {A sample session: Tuples of 4 involutions in $S_3$} 
 
gap$>$ \ g:=SymmetricGroup(3);;\\
gap$>$ \ t:=$[ (1,2), (1,2), (1,2), (1,2)]$;\\

\noindent
gap$>$ \  Braid(g,t);\\

\smallskip\noindent
Collecting 20 random tuples... done\\
 Cleaning done; 20 random tuples remaining\\

\medskip\noindent
Orbit 1:\\
Length=4\\
Generated subgroup size=6\\
Centralizer size=1\\
Remaining portion of structure constant=3 \\
Cleaning current orbit... done; 1 random tuples remaining\\

\medskip\noindent
Orbit 2:\\
Length=1\\
Generated subgroup size=2\\
Centralizer size=2\\
Remaining portion of structure constant=0\\
Cleaning current orbit... done; 0 random tuples remaining\\

\medskip\noindent
Summary: orbits of generating tuples\\

\smallskip\noindent
Orbit of Length 4\\
Inner Hurwitz curve genus = 0\\
Straight inner Hurwitz curve genus = 0\\

\part{Applications of the BRAID program}

\section {Brief explanation of the background on covers}\label{expl}

 Let $\P=\bC\cup\{\infty\}$ the Riemann sphere. 
A {\bf cover} of $\P$ (in the classical sense) is a compact
Riemann surface $X$ together with a non-constant analytic map
$f:X\to\P$ of finite degree. By Riemann's Existence Theorem,
$f$ can also be viewed as a morphism of complex algebraic curves.

Consider such a cover $f:X\to\P$ of degree $n$.
It has finitely many  branch points  $p_1,...,p_r\in\P$
(points whose preimage has cardinality less than $n$). Pick $p\in
\P\setminus\{p_1,...,p_r\}$, and choose loops $\gamma_i$ around
$p_i$ such that $\gamma_1,...,\gamma_r$ is a standard generating
system of the fundamental group $\Gamma:=\pi_1(
\P\setminus\{p_1,...,p_r\},p)$ (see \cite{Buch}, Thm. 4.27); in
particular, we have $\gamma_1\cdots\gamma_r=1$. Such a system
$\gamma_1,...,\gamma_r$ is called a homotopy basis of
$\P\setminus\{p_1,...,p_r\}$. The group $\Gamma$ acts on the fiber
$f^{-1}(p)$ by path lifting, inducing a transitive subgroup $G$ of
the symmetric group $S_n$ (determined by $f$ up to conjugacy in
$S_n$). It is called the {\bf monodromy group} of $f$. The images
of $\gamma_1,...,\gamma_r$ in $S_n$ form a 
tuple $\s=(\s_1,...,\s_r)$ generating $G$.
We say the  cover $f:X\to\P$ is {\bf of type} $\s$.
The genus $g$ of $X$
depends only on  $\s$, and is given by the {\bf Riemann-Hurwitz
formula}
$$2\ (n+g-1)\ \ \ \ = \ \ \ \sum_{i=1}^r\ \rm{Ind}(\s_i)
\leqno{(3)} $$
where the index $\rm{Ind}(\s_i)$ of a permutation in $S_n$ is
$n$ minus the number of orbits.

A tuple $\s=(\s_1,...,\s_r)$ of elements of
 $S_n$ arises in the above way
from a cover of degree $n$ if and only if $\s$ generates
a transitive subgroup $G$ and $\s_1\cdots \s_r=1$ and
$\s_i\ne1$ for all $i$. Call such a tuple {\bf admissible}.
The significance of  braid orbits comes from the following
fact (which follows from Nielsen's theorem).

\medskip\noindent
{\bf Theorem:} {\it Let $\s$ and $\s'$ be
admissible tuples generating the same subgroup $G$ of $S_n$.
Suppose $f:X\to\P$ is a cover  of type  $\s$.
Then $f$ is of type  $\s'$ if and only if   
the braid orbits of $\s$ and $\s'$ are conjugate under
$N_{S_n}(G)/G$.}

\medskip\noindent
Here $N_{S_n}(G)$ is the normalizer of $G$ in $S_n$.
The action of $N_{S_n}(G)/G$  on  braid orbits
comes from the fact that
if $\s$  generates  $G$ then 
Inn$(G)$ fixes the braid orbit of $\s$
 (see \cite{Buch}, Lemma 9.4). 

The next important fact is that 
the covers of type $\s$ form an {\it irreducible family}.
Here we use the term "family" in the non-technical sense:
Two covers are in the same irreducible family if they can be
continously deformed into each other (keeping the branch points
distinct). It turns out that the covers of type $\s$ are
parametrized (up to equivalence) by an irreducible variety,
the {\bf Hurwitz space} ${\cal H}(\s)$. 
This is made precise in the theory of Hurwitz spaces
(= moduli spaces for covers of $\P$), see  \cite{FV},
\cite{Buch},\cite{Vmod}. 

Two covers $f:X\to\P$ and $f':X'\to\P$ are called equivalent
(resp., weakly equivalent) if there is a homeomorphism
$h:X\to X'$ (resp., a homeomorphism $h:X\to X'$ and an 
analytic automorphism
$g$ of $\P$) such that $f=f'\circ h$ (resp., $g\circ f=f'\circ h$).
The automorphism group of $\P$ is PGL$_2(\bC)$ (group of
fractional linear transformations). It has a natural action on the
Hurwitz space ${\cal H}(\s)$. The quotient by this action is the
reduced Hurwitz space ${\cal H}^{\mbox{red}}(\s)$. 
It parametrizes the covers of type $\s$ up to weak equivalence.
Summarizing:

\medskip\noindent
{\bf Basic Fact:} {\it The covers of type $\s$ are
parametrized up to equivalence
(resp., up to weak equivalence) by an irreducible variety,
the {\bf Hurwitz space} ${\cal H}(\s)$
(resp., ${\cal H}^{\mbox{red}}(\s)$).
These varieties depend only on the braid orbit of $\s$.}

\medskip\noindent
A cover $f:X\to\P$ of type $\s$  is a Galois cover if and only if
$\s$ generates a regular subgroup $G$ of $S_n$. Pairs $(f,\mu)$,
where $f$  is a Galois cover of type $\s$ and $\mu: \mbox{Deck}(f)
\to G$ an isomorphism, are parametrized by the inner
Hurwitz space ${\cal H}_{\mbox{in}}(\s)$ (up to suitable equivalence).
This also is an irreducible variety.
Its quotient by PGL$_2(\bC)$ is the inner reduced Hurwitz space 
${\cal H}_{\mbox{in}}^{\mbox{red}}(\s)$. It is the inner
Hurwitz space that is of foremost importance for the Inverse Galois
Problem (see \cite{FV}). There is another version of it, the
straight inner
Hurwitz space $\tilde {\cal H}_{\mbox{in}}(\s)$ that  parametrizes
pairs $(f,\mu)$ together with an ordering of the branch points
of $f$. It also has a reduced version
 $\tilde {\cal H}_{\mbox{in}}^{\mbox{red}}(\s)$.

If $\s$ has length $r\le3$ then 
${\cal H}^{\mbox{red}}(\s)$ and ${\cal H}_{\mbox{in}}^{\mbox{red}}(\s)$
consist just of a single point.
If $r=4$ then these reduced Hurwitz spaces are curves. In the next
section we show how to compute their genus.

\section{The genus of the reduced Hurwitz curve in the case $r=4$}
             \label{genus}

In this section we look at the case $r=4$. 
The braid group $\BB_4=<Q_1,Q_2,Q_3>$ acts on Inn$(G)$-orbits 
of admissible 4-tuples from $G$ via its quotient $\overline \BB_4$ 
defined by the extra relations
$$Q_1Q_2Q_3^2Q_2Q_1\ \ = \ \ 1 \ \ = \ \ Q_1^2Q_3^{-2}$$
The structure of   $\overline \BB_4$
has been determined by
Thompson \cite{Th}. We denote the image of $Q_i$ in $\overline \BB_4$
by the same symbol, for simplicity. The elements $\gamma_0=Q_1Q_2$
and $\gamma_1=Q_1Q_2Q_1$ of $\overline \BB_4$ have order 3 and 2,
respectively. The elements $Q_1Q_3^{-1}$
and $(Q_1Q_2Q_3)^2$ generate a normal Klein 4-group
$\cal V$ in $\overline \BB_4$, and $\overline \BB_4/\cal V$ is the 
free product of the images of $<\gamma_0>$ and $<\gamma_1>$.

Fix an admissible 4-tuple $\s=(\s_1,...,\s_4)$, and
let $G\subset S_n$
be the  group generated by  $\s$.
Two 4-sets (unordered 4-tuples)
 of points of $\P$ are PGL$_2(\bC)$-conjugate
if and only if they have the same $j$-invariant (which can be any
complex number).  
The covers $f$ of type $\s$
 whose branch points have fixed $j$-invariant $\ne 0,1$
are parametrized, up to weak equivalence, by the set $F$
of $\cal V$-orbits
of  $N_{S_n}(G)$-orbits 
of  4-tuples in the braid orbit of $\s$. (Follows from the theory
outlined in section \ref{expl}, plus the fact that the stabilizer in
PGL$_2(\bC)$ of any 4-set with $j$-invariant $\ne 0,1$ is a 
Klein 4-group). From this one obtains an explicit description
of the 
{\bf Hurwitz curve} ${\cal H}^{\mbox{red}}(\s)$ parametrizing 
the covers of type $\s$ (up to weak equivalence).
It arises as
covering of $\P$ with branch points at $0,1,\infty$ whose general fiber
is in 1-1 correspondence with $F$.
The triple of permutations
 associated with this covering (by section \ref{expl}) is given
by the action on $F$ of $\gamma_0, \gamma_1$ and $\gamma_\infty:=Q_2$
(see \cite{BF}, Prop. 4.4 
and \cite{DF}, Prop. 6.5). From this we can compute the genus of 
${\cal H}^{\mbox{red}}(\s)$ by the Riemann-Hurwitz formula (2).
The case of the inner reduced Hurwitz curve
${\cal H}_{\mbox{in}}^{\mbox{red}}(\s)$ is analogous, with
$F$ replaced by the set 
of $\cal V$-orbits
of  Inn$(G)$-orbits
of  4-tuples in the braid orbit of $\s$.

\section {Indecomposable rational functions and 
primitive genus zero systems}

Here we are concerned with covers $f:X\to\P$ where $X$ has genus $0$.
Then we can identify $X$ with $\P$, so we consider 
covers $f:\P\to\P$. If such a cover has  degree $n$ then it is given
by a rational function of degree $n$, i.e., $f(x)=P(x)/Q(x)$ where
$P$ and $Q$ are complex polynomials with $n=\max(\deg(P),\deg(Q))$.
Then the monodromy group $G$ of $f$
is isomorphic (as a permutation group)
to the Galois group of the polynomial
$P(x)-tQ(x)$ over $\bC(t)$.
By the Riemann-Hurwitz
formula (2), genus $0$ covers correspond to the following
kind of tuples:

\begin{Definition}
A {\bf genus zero system} in $S_n$ is a tuple
$(\s_1,...,\s_r)$ generating a transitive subgroup $G$
of $S_n$ such that $\s_1\cdots\s_r=1$ and $\s_i\ne1$ (for all $i$)
and
  $$2\ (n-1)\ \ \ \ = \ \ \ \sum_{i=1}^r\ \rm{Ind}(\s_i)$$
It is called a {\bf primitive genus zero system} if $G$ is
primitive.
\end{Definition}

Thus by section \ref{expl}, irreducible families of
rational functions in $\bC(x)$
of degree $n$ with monodromy group $G\subset S_n$
correspond to $N_{S_n}(G)/G$-orbits of
 braid orbits of genus zero systems
generating $G$.
The family consists of indecomposable functions if and only if 
$G$  is primitive.
Here "indecomposable" means that the function is not the composition
$f_1(f_2(x))$ of two functions of degree $>1$.

There is a huge number of genus zero systems
that generate $S_n$ or $A_n$, too many
to be classified. The \lq general\rq\ rational function
has monodromy group $S_n$. Those functions with smaller monodromy group 
satisfy interesting identities and therefore it seems desirable to
have a complete classification of their irreducible families.
They correspond to the primitive genus zero systems that generate
a permutation group $G$ other than  $S_n$ or $A_n$. 
 The smallest case
is $G=L_3(2)$ (acting on 7 points). It has the most
 braid orbits of genus zero  systems. We discuss this example
in the following section.

\section{Example: Genus zero systems for the action of $G=L_3(2)$ on 7 points}
\label{L32}

The braid orbits of such tuples are listed in Table 1. 
We note  there is 
exactly one braid orbit $B_6$ of tuples of length 6, all the others
consist of shorter tuples. 

Replacing the last two entries of a tuple $\s$
 by their product is called
"Coalescing the tuple". Geometrically, this means that we merge
(or "coalesce") the
last two branch points of the associated cover. The family
corresponding
to the coalesced tuple $\s'$ lies in the boundary of the original family; 
in other words, the generic cover of type $\s'$
arises by specialization of the generic  cover of type $\s$. 

One checks that each of the orbits $\ne B_6$
in Table 1 contains a tuple that arises by a sequence of such
coalescing operations from a tuple of length 6. This means that there is
essentially only one
family of rational functions of degree 7 with
monodromy group $G=L_3(2)$. The generic function 
in this family has
6 branch points, and on the boundary we have functions with 3, 4 or 5
branch points. We can extract an explicit form of such a generic function
from  \cite{Ma}, Thm. 4.3:

\medskip\noindent

{\bf Generic function of degree 7 with
monodromy group} $L_3(2)$: \ \ $$f(x) \ \ = \ \ \frac{P(x)}{
x^2(x-c)(x^2-bx+b)}$$
where
$$P(x) \ \ = \ \ 
x^{7} - (a\,(c - 2) + 2\,b + c)\,x^{6} + ( \\
 - (b - 4)\,(c - 1)\,a^{2} + ((c - 2)\,b^{2} + (2\,c^{2} - 5\,c
 + 4)\,b - 2\,c^{2})\,a + b\,(2\,b\,c + 2\,c^{2} + b^{2}) \\
)x^{4}\mbox{} +
$$ $$
 ((2\,c^{2} - 1)\,(b - 4)\,a^{2} + (( - 2\,c^{2}
 + c + 2)\,b^{2} + (5\,c^{2} + 2\,c - 4)\,b - 4\,c^{2})\,a \\
\mbox{} - (c + 1)\,b^{3} - c\,(2\,c + 3)\,b^{2} + c^{2}\,b)x^{3}
\mbox{}
$$ $$ + ((c^{2} + 3\,c - 1)\,(4 - b)\,a^{2} \\
\mbox{} + ((3\,c - 2)\,b^{2} - 2\,(c^{2} + 4\,c - 2)\,b + 4\,c^{2
})\,a + b\,(b^{2} + 3\,b\,c - c^{2}))c\,x^{2} \\ $$ $$
\mbox{} + (2\,a\,b\,c - 8\,a\,c + a\,b - 4\,a - b^{2} + 2\,b\,c)
\,a\,c^{2}\,x - a^{2}\,(b - 4)\,c^{3}
$$

Replacing a function
$g(x)$ by $\alpha(g(\beta(x))$
with $\alpha ,\beta\in$ PGL$_2(\bC)$
doesn't change the monodromy group. So the functions we are
interested in are only determined up to coordinate change.
(Weak equivalence of covers, see above).

To illustrate the interplay between these functions and the group-theoretic
data in Table 1, we consider the specialization $b=0$.
The resulting  function $y=h(x)$ still has degree 7.
It has poles of order 4,2,1 at
$x=0,\infty,c$, respectively. Thus the corresponding tuple $\s$
contains an element of cycle type
$(4)(2)$ (corresponding to the branch point $y=\infty$).
Thus the monodromy group of $h(x)$ is still $L_3(2)$
(since it is  a transitive subgroup of $L_3(2)$ containing an element
of order 4).
The ramification index at a point $x=x_0$ not over $y=\infty$
equals one plus the multiplicity of  the zero $x=x_0$ of the
derivative $h'(x)$. Here we can replace $h'(x)$ by its numerator
(when it is written as a rational function in reduced form).
This numerator is a lengthy expression of degree 8 in $x$.
 But its discriminant with respect to $x$
factors nicely 
as $ 16777216\ c^{16}\,a^{9}$ times the cube of the following
expression (4) 
times the square of another (slightly longer) expression (5) that
we don't display here.
$$4\,a^{2}\,c^{4} + 8\,a\,c
^{4} + 4\,c^{4} - 4\,a^{2}\,c^{3} - 36\,c^{3}\,a + a^{3}\,c^{2}
 + 6\,a^{2}\,c^{2} \\
\mbox{} + 16\,c^{2}\,a - 2\,a^{3}\,c - 8\,c\,a^{2} + 2\,a^{3} +
16\,a^{2} \eqno{(4)}$$

The discriminant is non-zero, hence the above ramification indices
are all
 $\le2$. It follows that $\s$ consists of an element of order 4
and four involutions (by Riemann-Hurwitz). Thus $h(x)$ is the
generic function in the $(2A,2A,2A,2A,4A)$-family from Table 1.

Let's see how we can further specialize $h(x)$ by coalescing
two of the finite branch points. By Table 1, this leads to the
$(2A,2A,3A,4A)$- and the $(2A,2A,4A,4A)$-family. Both of those
have ramification indices $>1$ at certain points not over $y=\infty$.
Hence these specializations annihilate 
 the above discriminant. The two factors (4) and (5)
define genus zero curves in the $a,c$-plane (checked by 
\cite{Maple}).
This  corresponds nicely to the fact that the 
$(2A,2A,3A,4A)$- and the $(2A,2A,4A,4A)$-family are parametrized 
by  Hurwitz curves of genus zero (see Table 1).

Incidentally, \cite{Ma}, Thm. 4.2 gives another version of the
generic function in the $(2A,2A,2A,2A,4A)$-family.
(He doesn't consider our version).  One can
similarly specialize it to obtain two genus zero curves 
parametrizing the $(2A,2A,3A,4A)$- and the $(2A,2A,4A,4A)$-family.

\begin{table}[h] \label{l32}

\caption{Genus zero systems for the action of $G=L_3(2)$ on 7 points} 
\vskip 0.5cm 

\begin{center} 
\begin{tabular}{||c|c|c|c|c||} 
\hline \hline &&&&\\ 
classes $C_1,...,C_r$ & length of orbits & number of orbits & genus&straight  genus \\
\hline 
 $(2A,2A,2A,2A,2A,2A)$ & 1680 & 1  &&\\
\hline
 $(2A,2A,2A,2A,3A)$ & 216 & 1 & &\\
\hline 
 $(2A,2A,2A,2A,4A)$ & 192  & 1 & & \\
\hline
$(2A,2A,2A,7A)$ & 7 & 1 &0 &0 \\
\hline
$(2A,2A,2A,7B)$ & 7 & 1 &0 &0\\
\hline
$(2A,2A,3A,3A)$ & 30 & 1 &0 &2\\
\hline
$(2A,2A,3A,4A)$ & 24 & 1 &0 &1\\
\hline
$(2A,2A,4A,4A)$ & 24 & 1 &0 &1\\
\hline
$(2A,3A,7A)$ & 1 & 1 & &\\
\hline
$(2A,3A,7B)$ & 1 & 1 & &\\
\hline
$(2A,4A,7A)$ & 1 & 1 & & \\
\hline
$(2A,4A,7B)$ & 1 & 1 &  &\\
\hline
$(3A,3A,4A)$ & 1 & 4&  &\\
\hline
$(3A,4A,4A)$ & 1 & 2&  &\\
\hline
$(4A,4A,4A)$ & 1 & 4&  &\\

\hline \hline

\end{tabular} 
\end{center} 
\end{table} 

\section{Primitive genus zero covers 
 branched at $\ge5$ points} \label{five}

Each primitive permutation group has a characteristic subgroup
$F^*(G)$ (called the generalized Fitting subgroup)
which is  the 
 direct product of isomorphic simple groups.
Frohardt, Guralnick and Magaard \cite{FGM2} determine all
primitive genus zero systems generating a group $G\subset S_n$
with $F^*(G)$ not abelian and not 
a direct product of alternating groups. 
The resulting list is finite,
but too long to be shown in tabular form. However,
there are only a few cases with $r\ge5$ (i.e., where the corresponding
covers  are branched
at $5$ or more points). We list these in Table 2 and note that for each 
choice
of $C_1,...,C_r$ there is exactly one associated braid orbit
(i.e., exactly one irreducible family of genus zero covers).

The table was produced as follows. A series of reductions shows
that the permutation degree of such a system is at most $1000$.
It remains to search the GAP library of 
primitive permutation groups of degree $\le1000$.
For each such group $G$  that
satisfies our hypothesis, we find all collections of
 conjugacy classes $C_1,...,C_r$ 
 that satisfy the Riemann-Hurwitz
formula (for $g=0$). For each such collection, we apply the
BRAID program to find all braid orbits of associated tuples.

\begin{table}[h] \label{rk 5} 

\caption{Braid orbits of primitive genus zero systems of length $\ge5$
in almost simple groups $\ne A_n,S_n$}
\vskip 0.5cm 

\begin{center} 
\begin{tabular}{||c|c|c|c||} 

\hline \hline &&& \\ 
$G$ & degree & classes $C_1,...,C_r$ & orbit length \\
\hline
$L_4(3)$ & 40 & $(2A,2B,2B,2C,2C)$ & 320  \\
\hline
$S_6(2)$ & 36 & $(2A,2B,2B,2B.3B)$ &  4  \\
\hline 
$L_5(2)$ & 31 & $(2B,2B,2B,2B,2B)$ & 31744 \\
\hline
& 31 & $(2A,2A,2B,2B,3B)$ & 528 \\
\hline
$S_6(2)$ & 28 & $(2A,2A,2A,3B,4A)$ &  4 \\
\hline 
& 28 & $(2A,2C,2C,2C,3B)$ &  54 \\
\hline 
& 28 & $(2A,2D,2D,2D,2D)$ &  3584 \\
\hline
$M_{24}$ & 24 & $(2A,2A,2A,2A,4B)$ &  72000 \\
\hline
$M_{23}$ & 23 & $(2A,2A,2A,2A,3A)$ &  21456 \\
\hline
$M_{22}$ & 22 & $(2A,2A,2A,2B,2C)$ & 660 \\
\hline
& 22 & $(2A,2A,2B,2B,3A)$ & 600 \\
\hline 
$L_3(4)$ & 21 & $(2A,2A,2A,2A,2A)$ & 252 \\
\hline
$L_3(4).3.2_2$ & 21 & $(2B,2B,2B,2B,3A)$ & 1824 \\
\hline 
 & 21 & $(2A,2A,2B,2B,3B)$ & 264 \\
\hline
$L_3(3)$ & 13 & $(2A,2A,2A,2A,2A,2A)$ & 32760 \\
\hline 
& 13 & $(2A,2A,2A,2A,3B)$ & 1944 \\
\hline
& 13 & $(2A,2A,2A,2A,4A)$ & 2016 \\
\hline
& 13 & $(2A,2A,2A,2A,6A)$ & 2160 \\
\hline 
& 13 & $(2A,2A,2A,3A,3A)$ & 120 \\
\hline 
$M_{12}$ & 12 & $(2A,2A,2A,2A,2B)$ & 2048 \\
\hline
& 12 & $(2A,2A,2A,2A,3A)$ & 2784 \\
\hline 
& 12 & $(2A,2A,2A,2A,4B)$ & 7296 \\
\hline
$M_{11}$ & 12 & $(2A,2A,2A,2A,3A)$ & 2376 \\
\hline
$L_2(11)$ & 11 & $(2A,2A,2A,2A,2A)$ & 704 \\
\hline
$L_3(2)$ & 7 & $(2A,2A,2A,2A,2A,2A)$ & 1680 \\
\hline
& 7 & $(2A,2A,2A,2A,3A)$ & 216 \\
\hline 
& 7 & $(2A,2A,2A,2A,4A)$ & 192 \\
\hline 
\hline 
\end{tabular} 
\end{center} 
\end{table}


\begin{thebibliography}{99}

\bibitem[BF] {BF} {\sc P. Bailey and M.~Fried},
Hurwitz monodromy, spin separation and higher levels of a
modular tower,
Proceedings of Symposia in Pure Math. {\bf 70} (2002),
79--220.

\bibitem[Br] {Br} {\sc  Th. Breuer},
Characters and automorphism groups of compact Riemann surfaces,
London Math. Soc. Lect. Notes {\bf 280}, Cambridge Univ. Press 2000.

\bibitem[DF] {DF} {\sc P. Debes and M.~Fried},
Integral specialization of families of rational functions,
Pacific J. Math. {\bf 190} (1999), 45--85.


\bibitem[De] {Dett} {\sc  M. Dettweiler},
Kurven auf Hurwitzr\"aumen und ihre Anwendungen in der Galoistheorie, 
Dissertation, Erlangen, 1999

\bibitem[FG] {FG} {\sc M.~Fried and R. Guralnick},
On uniformization of generic curves of genus $g<6$ by radicals,
unpublished manuscript.

\bibitem[FV] {FV} {\sc M.~Fried and H. V\"olklein},
 The inverse Galois problem and
 rational points on moduli spaces,   {\sl Math. Annalen}
 {\bf  290} (1991), 771--800.

\bibitem[FGM1] {FGM} {\sc D. Frohardt, R. Guralnick and K. Magaard},
Genus zero actions of groups of Lie rank 1, Proc. Symp. Pure Math.
{\bf 70} (2002), 449--483. 

\bibitem[FGM2] {FGM2} {\sc D. Frohardt, R. Guralnick and K. Magaard},
The primitive genus zero systems involving non alternating,
non abelian simple groups, preprint

\bibitem[FM] {FM} {\sc D. Frohardt and K. Magaard},
Composition factors of monodromy groups,
Annals of Math. {\bf 154} (2001),1-19.

\bibitem[GAP4] {GAP4} {\sc The  GAP~Group},
\emph{ GAP -- Groups, Algorithms, and Programming}, Version 4.2; 2000.
\verb+(http://www.gap-system.org)+

\bibitem[Gra] {Gra}  {\sc L. Granboulan},
Construction d'une extension r\'eguli\`ere de $\bQ(t)$ de groupe
de Galois $M_{24}$, Exp. Math. {\bf 5} (1996), 3--14.

\bibitem[GM] {GM} {\sc R. Guralnick and K. Magaard},
On the minimal degree of a primitive permutation group,
J. Algebra {\bf  207} (1998), 127--145.

\bibitem[GN] {GN} {\sc R. Guralnick and M. Neubauer},
Monodromy groups and branched coverings: The generic case,
Contemp. Math. {\bf 186} (1995), 325--352.

\bibitem[GS] {GS} {\sc R. Guralnick and J. Shareshian},
Alternating and Symmetric Groups as Monodromy
Groups of Curves I, preprint.

\bibitem[MSSV] {kyoto}  {\sc K. Magaard, S. Shpectorov and H. V\"olklein},
 The locus of curves with prescribed automorphism group
RIMS Publication series {\bf 1267} (2002), 112--141
 (Communications in Arithmetic Fundamental Groups,
Proceedings of the RIMS workshop held at Kyoto University Oct. 01)

\bibitem[MV] {MV}  {\sc K. Magaard and H. V\"olklein},
The monodromy group of a function on a general curve, submitted
(see also math.AG/0304130).

\bibitem[Ma]{Ma}
 {\sc G. Malle}, 
Multi-parameter polynomials with given Galois group,
J. Symb. Comp. {\bf 30} (2000), 717--731.

\bibitem[MM]{MM}
 {\sc G. Malle and B. H. Matzat}, Inverse Galois Theory,
Springer, Berlin-Heidelberg-New York 1999.

\bibitem [Maple] {Maple} {\sc Maple 6}, Waterloo Maple Inc.,  2000.

\bibitem[Ne1] {Ne1} {\sc M. Neubauer},
On primitive monodromy groups  of genus 0 and 1,
Comm. Alg. {\bf 21} (1993), 711--746.


\bibitem[Ne2] {Ne2} {\sc M. Neubauer},
On monodromy groups  of fixed genus,
J. Algebra {\bf 153}(1992), 215--261.

\bibitem[Prz] {Pr} {\sc B. Przywara},
Braid operation software package 2.0 (1998),
available at http://www.iwr.uni-heidelberg.de/ftp/pub/ho

\bibitem[Th] {Th}  {\sc J. Thompson },
Note on H(4), Comm. Alg. {\bf 22}(1994), 5683--5687.


\bibitem[ThV] {ThV}  {\sc J. Thompson and H. V\"olklein}, 
Symplectic groups as Galois groups, {\em J. Group Theory}
{\bf 1} (1998), 1--58.

\bibitem[V1] {Buch} {\sc  H.~V\"olklein},
 Groups as Galois Groups --
an Introduction, Cambr. Studies in Adv.
Math. 53, Cambridge Univ. Press 1996.

\bibitem[V2] {Vmod} {\sc H. V\"olklein},
 Moduli spaces for covers of the Riemann sphere,
 Israel J. Math.  {\bf 85} (1994), 407--430.

\end{thebibliography}
\end{document}